\documentclass[11pt,a4paper]{article}
\usepackage[reqno,sumlimits]{amsmath}
\usepackage{latexsym}
\usepackage{graphics}
\usepackage{amsfonts} 
\usepackage{amssymb}   
\usepackage{theorem}
\usepackage{epic,eepic}
\usepackage{epsfig}
\usepackage{amscd}
\newtheorem{theorem}{Theorem}[section]
\newtheorem{proposition}[theorem]{Proposition}

\newtheorem{lemma}[theorem]{Lemma}

{\theoremstyle{plain}
{\theorembodyfont{\normalfont\rmfamily}
\newtheorem{definition}[theorem]{Definition} 

}}
\newcommand{\cqd}{\hfill $\Box$\medskip}
\newcommand{\proof}{\noindent{\bf Proof: }}
\newcommand{\abs}[1]{\left| #1\right|}

\newcommand{\id} {\mbox{\rm Id}}
\newcommand{\im} {\mbox{\rm Im\,}}

\newcommand{\tr} {\mbox{\rm tr\,}}

\newcommand{\A}{{\cal A}}

\newcommand{\DD}{{\cal D}}

\newcommand{\II}{{\cal I}}
\newcommand{\JJ}{{\cal J}}

\newcommand{\LL}{{\cal L}}

\newcommand{\OO}{{\cal O}}

\newcommand{\RR}{{\cal R}}
\newcommand{\TT}{{\cal T}}
\newcommand{\UU}{{\cal U}}

\newcommand{\Cc}{{\mathbb{C}}}

\newcommand{\Ee}{{\mathbb{E}}}

\newcommand{\Ii}{{\mathbb{I}}}

\newcommand{\Pp}{{\mathbb{P}}}
\newcommand{\Qq}{{\mathbb{Q}}}
\newcommand{\Rr}{{\mathbb{R}}}

\newcommand{\Tt}{{\mathbb{T}}}

\newcommand{\Vv}{{\mathbb{V}}}

\newcommand{\Zz}{{\mathbb{Z}}}
\begin{document}
\begin{titlepage}

\title
{\bf Renormalisation scheme for vector fields on $\Tt^2$ with a diophantine frequency}

\author{Jo\~ao Lopes Dias\thanks{Email: j.lopes-dias@damtp.cam.ac.uk}\\   
{\small Department of Applied Mathematics and Theoretical Physics,} \\
{\small University of Cambridge} \\ 
{\small Silver Street, Cambridge  CB3 9EW, England}}
\date{April, 2001}

\maketitle
\begin{abstract}
We construct a rigorous renormalisation scheme for analytic vector
fields on $\Tt^2$ of Poincar\'e type.
We show that iterating this procedure there is convergence to a limit
set with a ``Gauss map'' dynamics on it, related to the continued fraction
expansion of the slope of the frequencies.
This is valid for diophantine frequency vectors.
\end{abstract}

\end{titlepage}

%%%%%%%%%%%%%%%%%%%%%%%%%%%%%%%%%%%%%%%%%%%%%%%%%%%%
\section{Introduction}

A renormalisation operator $\RR$ is defined for analytic vector fields on the
torus of dimension two, associated with a frequency vector.
This approach is based in \cite{Koch} for Hamiltonian functions, and in
\cite{jld} for flows on $\Tt^d=\Rr^d/\Zz^d$, $d\geq2$.
The vector fields considered here are of Poincar\'e type, i.e. there is
a classification by a unique winding ratio on which $\RR$ acts as the Gauss
map.
The slope $\alpha$ of the winding ratio 
is mapped by the
shift of its continued fraction expansion.
All vector fields corresponding to 
non-zero $\alpha$ are
renormalisable, and those with irrational $\alpha$ are
infinitely renormalisable.
We extend a result in \cite{jld} based on \cite{Koch} for diophantine winding ratios.

The renormalisation group technique has been used for some time
applied to various problems in dynamical systems (e.g. \cite{Feigenbaum,MacKayThesis,Stark}).
Roughly, it analyses systems on longer time
scales and smaller spacial scales, generating a new system.
By iterating this we may get convergence to a limiting behaviour,
preferably simple, e.g. a fixed point.
In such a case, we get self-similarity,
which can be non-trivial depending on the problem.

We show here that an orbit by $\RR$ of a
constant vector field with diophantine winding ratio, attracts all the nearby
orbits in the same homotopy class.
This can be applied to other flows
on domains with an extra vertical dimension.
Also, it might be used for two-degrees of freedom
Hamiltonian systems, in order to extend \cite{Koch}.
The allowed set of frequencies is smaller than the one obtained by
KAM theory.
Nevertheless, this is not a real disadvantage of our method since it
is of equal full Lebesgue measure.

Our procedure can be related to the one followed by MacKay
\cite{MacKay2}, and Chandre and Moussa \cite{Chandre} for an
approximate renormalisation scheme in the
framework of the
breakup of invariant tori in Hamiltonian systems.
The renormalisation studied by those authors, besides truncating the
Fourier modes of the Hamiltonians, does not consider perturbations
depending on the action variable. 
That would have implied a change in
the Hamiltonian vector fields, namely the angle time derivative.
Thus, it does not have the unstable
direction that appears in our systems corresponding to perturbations
non-collinear to $\omega$.

\medskip

In Section \ref{change of basis for generic freq},
we include the main idea for the construction of the linear part of
the renormalisation, along with the main properties of the continued
fraction expansion of an irrational number.
Then, in Section \ref{Space of Analytic Vector Fields} we introduce the 
space of vector fields of our interest.
We define a one-step operator $\RR$ in Section \ref{One-step
Renormalisation Operator} and show one of its statements in Section \ref{section sets
I+ and Ik}.
It follows Section \ref{Diophantine Numbers} with more known results of number theory, this time on
the diophantine condition.
The main convergence theorem for diophantine vectors is presented in Section \ref{The Limit Set of the
Renormalisation}, along with a proof in Section \ref{proof of theorem generic frequency d=2}.

%%%%%%%%%%%%%%%%%%%%%%%%%%%%%%%%%%%%%%%%%%%%%%%%%%%%%%%%%%%%%%%%%%%%%%%%%%%%

\section{Change of Basis for a Generic Frequency and the Continued
Fraction Expansion} \label{change of basis for generic freq}

It was shown in \cite{jld} (using results from \cite{Koch}) that for a
quadratic irrational slope of a two-dimensional vector $\omega$, the
renormalisation operator therein constructed is well-defined since its linear part
exists.
More precisely, we can find a hyperbolic matrix $T$ in $GL(2,\Zz)$ for
which $\omega$ is an eigenvector.
For other vectors one needs to construct a
different operator.
In fact, for each iteration of the renormalisation we will
consider a different transformation, determined by the
continued fraction expansion of $\alpha$.
%The limit set depends on its arithmetic properties, as will be seen.

\medskip

Consider the constant vector field $\omega$ on the universal cover $\Rr^2$,
with slope $\alpha=\omega_2/\omega_1$, assuming that $\omega_1\not=0$.
The main idea is to perform a change of basis enlarging the region of
$\Rr^2$ around the orbits, and a scaling of time.

In the case $\alpha>1$, we choose the new basis $\{(1,1),(0,1)\}$, that
decreases (``$D$'') the slope of the transformed vector: $\alpha'=\alpha-1$. 
On the other hand, if $0<\alpha<1$, one simply swaps (``$S$'') the coordinates,
i.e. we choose $\{(0,1),(1,0)\}$ thus inverting the slope: $\alpha'=1/\alpha$.
The change of basis matrices are then
$D=\left[\begin{smallmatrix}1&0\\1&1\end{smallmatrix}\right]$ and
$S=\left[\begin{smallmatrix}0&1\\1&0\end{smallmatrix}\right]$,
respectively.
% (see Figures \ref{action of D} and \ref{action of S}).
Note also that
$D^a=\left[\begin{smallmatrix}1&0\\a&1\end{smallmatrix}\right]$,
$a\in\Zz$ and $S^{-1}=S$.
Finally, if $\alpha<0$, we recover one of the previous cases using $V=\left[\begin{smallmatrix}-1&0\\0&1\end{smallmatrix}\right]=V^{-1}$.

Recall that the matrices $D,S$ and $V$ generate $GL(2,\Zz)$. 
This group acts on $\bar \Rr=\Rr\cup\{\pm\infty\}$ as $\alpha\mapsto (c+d\alpha)/(a+b \alpha)$ with $\left[\begin{smallmatrix}a&b\\c&d\end{smallmatrix}\right]\in GL(2,\Zz)$.
Each of the matrices $D,S$ and $V$ correspond to a transformation on
$\bar\Rr$ (in fact their inverses since those are the direct transformations to
the new basis): $D$ and $S$ as above and $V\colon\alpha'=-\alpha$.

The {\it continued fraction expansion}\index{continued fraction
expansion} of $\alpha\in\Rr$, is given by
$$
\alpha=a_0+
\cfrac1{a_1+\cfrac1{a_2+\dots}},
$$
or simply $\alpha=[a_0,a_1,a_2,\dots]$, with $a_0\in\Zz$ and
$a_n\in\Zz^+$, $n\geq1$.
It is an infinite sequence if and only if $\alpha$ is irrational
\cite{Lang}.
We assume in the following that that is case and $\alpha>0$.
In these conditions, the continued fraction expansion is a one-to-one map.

Let $[x]=\max\{k\in\Zz\colon k\leq x\}$ 
be the integer part of $x\in\Rr$
and 
$\{x\}=x-[x]$ 
the fractional one.
The Gauss map\index{Gauss map} is given by 
$$
G\colon x\mapsto \left\{\frac1x\right\},\quad x>0.
$$ 
Writing $x_n=G(x_{n-1})$, $n\geq1$, with $x_0=\{\alpha\}$, one obtains
the above coefficients $a_n$
from the recurrence relation 
$a_0=[\alpha]$, $a_n=[1/x_{n-1}]$, $n\geq1$.
Note that $\alpha_n=a_n+x_n$, $n\geq0$, where
$\alpha_n=[a_n,a_{n+1},\dots]$, hence $\alpha_{n+1}=1/\{\alpha_n\}=1/x_n$.
Define also 
\begin{equation}\label{definition of beta_n}
\beta_n=\prod_{i=0}^n x_i \,= \prod_{i=0}^n \frac1{\alpha_{i+1}}.
\end{equation}
If $x_k\geq \gamma^{-1}$ for some $k=0,\dots, n-1$, where
$\gamma=\frac{1+\sqrt5}2$, then, letting $m=x_k^{-1}-x_{k+1}\geq 1$,
$x_kx_{k+1}=1-mx_k\leq 1-x_k \leq 1-\gamma^{-1}=\gamma^{-2}$.
On the other hand, if $x_k\leq\gamma^{-1}$, then
$x_{k+1}\geq\gamma^{-1}$, and the product is again less or equal to $\gamma^{-2}$.
Hence, 
\begin{equation}\label{bound on beta_n}
\beta_n\leq \gamma^{-n},
\quad
\frac{\beta_n}{\beta_{j-1}}\leq\gamma^{-(n-j)},
\quad
0\leq j\leq n.
\end{equation}

The quadratic irrationals\index{quadratic irrational} (roots of a quadratic polynomial over $\Qq$)
have eventually periodic continued
fraction expansion \cite{Lang},
e.g. $\gamma=[1,1,1,\dots]$ and $\sqrt{2}=[1,2,2,\dots]$.

The action of $D^{a_0}$ on $\alpha$ is the elimination of the
first coefficient $a_0$ of the continued fraction expansion,
$[a_0,a_1,a_2,\dots]\mapsto[0,a_1,a_2,\dots]$.
The transformation $S$ acts on sequences of the form
$[0,a_1,a_2,\dots]$ giving $ [a_1,a_2,\dots]$.
The composition is simply the shift of the coefficients.

We are looking for a specific type of linear coordinate changes of the
torus.
That is, it has to correspond to a matrix capable of improving the
analyticity, as will be seen in the following sections (see also \cite{jld}).
So, it has to be a hyperbolic matrix in $GL(2,\Zz)$, with an unstable direction close
enough to the subspace spanned
by $\omega$.
It turns out that we can choose the above shift of coefficients.
The complete linear transformation is then
$$
T_a=D^{a}S=\left[\begin{smallmatrix}0&1\\1&a\end{smallmatrix}\right]\in
GL(2,\Zz),
$$ 
where $a=[\alpha]$.
If $a\geq1$, this is a hyperbolic matrix with two real eigenvalues:
one inside and the other outside the unit circle.
More specifically,
$$
\lambda=\frac{a+\sqrt{a^2+4}}{2}\geq\gamma
$$ 
and $-1/\lambda$ are the eigenvalues.
The corresponding unstable and stable eigenvectors are $(1,\lambda)$ and
 $(1,-1/\lambda)$.

Now, it is useful to recall some properties of the continued fraction
expansion, since we have associated it with the change of basis to be
used in the construction of the renormalisation.
We can write the above continued fraction expansion procedure in
terms of the ``convergent'' matrices\index{convergent matrices}
$$
P_n=\left[ 
\begin{matrix} 
q_{n-1} & p_{n-1} \\
q_{n} & p_{n}
\end{matrix}
\right],\quad n\geq0,\quad\text{ and }
P_{-1}=\left[
\begin{matrix}
1 & 0\\
0 & 1
\end{matrix}
\right],
$$
using the recurrence formula 
$P_{n+1}=T_{n+1}P_n$,
where $T_n=T_{a_n}$, and
\begin{equation}\label{recurrence formulae q_n, p_n}
\begin{array}{l}
p_{n+1}=a_{n+1}p_n+p_{n-1}, \\
q_{n+1}=a_{n+1}q_n+q_{n-1}, \\
\beta_{n-2}=a_n\beta_{n-1}+\beta_n,
\end{array}
\end{equation}
assuming $\beta_{-1}=1$.
The last equality is simply derived from rewriting
$\beta_{n-2}=\alpha_{n}\beta_{n-1}$.
Notice that $T^*_{n+1}=(P_n^*)^{-1}P_{n+1}^*$ and 
$P_n=T_n\cdots T_0$, where $T^*$ stands for the transpose matrix of $T$.

The continued fraction expansion can be obtained in only one step by a
change of basis to $\{(q_{n-1},p_{n-1}),(q_n,p_n)\}$,
with the corresponding matrix being $P^*_n$. 
Thus, 
$(1,\alpha_{n+1}) = (q_{n-1}-\alpha q_n) P^*_n(1,\alpha)$. 
So, 
\begin{equation}\label{formula for alpha}
\alpha=\frac{p_{n-1}+p_n\alpha_{n+1}}{q_{n-1}+q_n\alpha_{n+1}},
\end{equation}
which implies $x_n=-(q_n\alpha-p_n)/(q_{n-1}-p_{n-1})$ and, by (\ref{definition of beta_n}),
$$
\beta_n=(-1)^{n}(\alpha q_n-p_n).
$$
Again, from (\ref{formula for alpha}) using $(n+1)$ instead of $n$,
$|\beta_n| = 1/(q_{n+1}+q_nx_{n+1})$. 
Hence,
\begin{equation}\label{upper and lower bounds on beta_n}
\frac{1}{2q_{n+1}} < \beta_n < \frac1{q_{n+1}}.
\end{equation}

It is easy to see that $\det P_n=(-1)^{n+1}$ and each line of the matrices
$P_n\in GL(2,\Zz)$ gives the convergents $p_n/q_n=[a_0,\dots,a_n]\to\alpha$, as
$n\to\infty$. 
This sequence verifies the best diophantine approximation
property\index{diophantine!approximation},
$|p_n-\alpha q_n|<|p-\alpha q|$, for any integer $p$ and $0<q<q_n$ (\cite{Hardy} - Theorem 182).
Also, if $|\alpha-p/q|<1/2q^2$, then $p/q$ is a convergent of $\alpha$
(\cite{Hardy} -- Theorem 184).
Using the fact that $|\det P_n|=1$, it is clear that the
convergents verify $|p_n/q_n-p_{n+1}/q_{n+1}|=1/q_{n}q_{n+1}$.
Then,
\begin{equation}\label{bounds on (alpha-p/q) and beta_n}
\abs{\alpha-\frac{p_n}{q_n}} <\frac1{q_{n}q_{n+1}}.
\end{equation}

Denote the product of the first $(n+1)$ continued fraction coefficients of
an irrational $\alpha$ by $A_n=\prod_{i=0}^{n}a_i$. 
Rewriting (\ref{recurrence formulae q_n, p_n}) one obtains
$q_{n+1}/q_n= a_{n+1} + 1/ (q_n/q_{n-1})$. By induction it can be
proved that the following is true:
$$%\begin{equation}\label{bound A and P on q_n}
A_n\leq q_n \leq A_n\, \prod_{i=1}^{n}(1+1/a_ia_{i-1}).
%\end{equation}
$$
Similarly, we define by $\tilde A_n$ the product of the $(n+1)$ first numbers
$\alpha_i = [a_i, a_{i+2}, \dots]$:
$$
\tilde A_n=\prod_{i=0}^n\alpha_i.
$$
This is related trivially to $\beta_n$ by the formula $\tilde
A_{n+1}=\alpha_0/\beta_n$.
So, from (\ref{recurrence formulae q_n, p_n}), $\tilde
A_{n-1}^{-1}=a_0\tilde A_n^{-1}+\tilde A_{n+1}^{-1}$, and, from
(\ref{upper and lower bounds on beta_n}),
\begin{equation}\label{upper and lower bounds on tilde A(n)}
\alpha_0 q_n < \tilde A_n < 2\alpha_0 q_n.
\end{equation}
Also, using (\ref{bound on beta_n}), 
\begin{equation}\label{bound on tilde A(n)}
\tilde A_n \geq \alpha_0 \gamma^{n-1},
\quad
\frac{\tilde A_n}{\tilde A_{j-1}}\geq \gamma^{n-j},
\quad
0\leq j\leq n.
\end{equation}

%%%%%%%%%%%%%%%%%%%%%%%%%%%%%%%%%%%%%%%%%%%%%%%%%%%%%%%%%%%%%%%%%%%%%%%%%%%%

\section{Space of Poincar\'e Analytic Vector Fields}\label{Space of Analytic Vector Fields}

Let $\phi_t$ be the flow generated
by a continuous vector field
$X$ on $\Tt^2$,
$\dot\theta=X(\theta)$, $\theta\in \Tt^2$, and
$\Phi_t$ its lift to $\Rr^2$.
For some norm $\|\cdot\|$ in $\Rr^2$, we say that 
$w_X(\theta_0)=\lim_{t\to\infty} \Phi_t(\theta_0)/
\|\Phi_t(\theta_0)\|$
is the {\it winding ratio} of $X$ for the orbit of
$\theta_0\in\Tt^2$, if the limit exists and $\lim_{t\to\infty}\|\Phi_t(\theta_0)\|=\infty$.
Otherwise, if $\Phi_t(\theta_0)$ is bounded we put $w_X(\theta_0)=0$.
On the other hand, if the limit does not exist or if
$\|\Phi_t(\theta_0)\|$ is unbounded but does not
tend to infinity, we do not define the winding ratio.

In the two-dimensional case, the set of winding ratios of vector
fields $X$ on $\Tt^2$ is always a subset of
$\{w,0,-w\}$, for some normalised vector $w$ (cf. e.g. \cite{Anosov,BGKM}).
Such a set is called {\it winding set}, and will be denoted by $w(X)$.

For any vector field $X$, the existence of an equilibrium is
equivalent to $0\in w(X)$, because it
corresponds to a bounded orbit.
If the slope of $w\in w(X)$ is an irrational number $\alpha$, 
then $w(X)=\{w\}$ (see \cite{BGKM} and references therein).
Fixed-point-free flows on $\Tt^2$ with only one winding ratio are called {\it Poincar\'e flows}.

A closed cross section or {\it transversal} to a vector field on a
surface is a simple closed curve, such that the vector field is nowhere
tangent to the curve, and intersects all the orbits.
A flow on $\Tt^2$ has a transversal if and only if it is a
Poincar\'e flow.
In this case, considering the return map on the transversal, it is possible to reduce the flow to a
diffeomorphism of the circle. 
Therefore, there is an equivalence between these systems and the
results are directly related.
If $w$ is rational (more precisely its slope), there is a periodic orbit of
winding ratio $w$ and every orbit is asymptotic to such.
When $w$ is irrational and the flow is $C^2$, then all the orbits are
dense on $\Tt^2$, and the flow is topologically equivalent to a
constant flow with the same winding ratio \cite{Denjoy,MacKay3}.

\medskip

We briefly recall the space of vector fields used in \cite{jld}
that we will also study here for the two-dimensional case.

Let $r>0$ and the lift of $\Tt^2$ to a complex
neighbourhood of $\Rr^2$:
$$
\DD(r)=\left\{\theta\in\Cc^2 \colon \|\im\theta\| <\frac r{2\pi}\right\},
$$
where $\|\cdot\|$ is the $\ell_1$-norm on $\Cc^2$.
That is,
$\|z\|=|z_1|+|z_2|$, with 
$|\cdot|$ the usual norm on $\Cc$.
We also denote the inner product as
$z\cdot z'=z_1z_1'+z_2z_2'$.

We will be dealing with analytic vector fields 
of the form
$X=\omega+f$, where $\omega\in\Rr^2$ 
and with analytic functions
$f\colon\DD(r)\to\Cc^2$, $2\pi$-periodic in each
variable $\theta_i$, represented in
Fourier series as
$$
f(\theta)=\sum\limits_{k\in\Zz^2}f_ke^{2\pi i k\cdot \theta}, \quad
f_k \in \Cc^2.
$$
More specifically, $f$ belongs to one of the
spaces $(\A(r),\|\cdot\|_r)$ and
$(\A'(r),\|\cdot\|'_r)$ whose
elements are
such that the respective norms
$$
\|f\|_r=\sum\limits_{k\in\Zz^2}\|f_k\|e^{r\|k\|}
\quad
\text{ and }
\quad
\|f\|'_r=\sum\limits_{k\in\Zz^2}\left(1+2\pi\|k\|\right)\|f_k\|e^{r\|k\|}
$$
are finite.
Both are Banach spaces.

Some of the above vector fields $X=\omega+f$ generate Poincar\'e flows
on $\Tt^2$.
Note that, if $\|f\|'_r < \|\omega\|$ there is no equilibria.
We will be interested in Poincar\'e flows, especially with
``irrational'' winding ratio.

%%%%%%%%%%%%%%%%%%%%%%%%%%%%%%%%%%%%%%%%%%%%%%%%%%%%%%%%%%%%%%%%%%%%%%%%%%%%

\section{One-step Renormalisation Operator}\label{One-step Renormalisation Operator}

In this section fix a vector $\omega=(\omega_1,\omega_2)\in\Rr^2\setminus\{0\}$ with
$\omega_1\not=0$ and
slope $\alpha=\omega_2/\omega_1\geq1$, and denote $a=[\alpha]$.
As in \cite{Koch}, we separate the terms of a vector field $X$ into
the following:

\begin{definition}\label{far from resonance terms}
For $\sigma>0$ and $\psi\in\Cc^2$, we define the {\it far
from resonance terms}\index{far
from resonance terms} with respect to $\psi$ to be the ones whose indices are in
$$
I_\sigma^-(\psi)=\left\{ k\in\Zz^2\colon |\psi\cdot k| > 
\sigma\|k\|\right\}.
$$
Similarly, the {\it resonant terms}\index{resonant terms} are
$I_\sigma^+(\psi)=\Zz^2\setminus I_\sigma^-(\psi)$.
The respective projections
$\Ii_\sigma^+(\psi)$ and $\Ii_\sigma^-(\psi)$
on the space of vector fields are defined as
$$
[\Ii_\sigma^\pm(\psi)] X(\theta)=\sum\limits_{k\in I_\sigma^\pm(\psi)}X_ke^{2\pi ik\cdot\theta}.
$$ 
\end{definition}

Let
$\omega'=\{\alpha\}^{-1}T^{-1}_{a}\omega$, with slope
$\alpha'=\{\alpha\}^{-1}$ that corresponds to the shift of the
continued fraction expansion of $\alpha$,
so that $\|\omega'\|/\|\omega\|=(1+\alpha')/(1+\alpha)$.
Denote also by
$\hat\omega'=\omega'/(\omega'\cdot\omega')$ the
normalised vector in the subspace spanned by $\omega'$, i.e.
$\hat\omega'\cdot\omega'=1$.
We denote by $\Ee(X)$ the spatial average of a vector field,
$\Ee(X)=\int_{\Tt^2}X(\theta)\,d\theta$, where $d\theta$ is the
normalised Lebesgue measure.

\begin{proposition}\label{proposition L}
There exists $0<\kappa<1$ and $\sigma>0$ such that,
if $0<\rho'<\rho$ with $\kappa\rho<\rho'$,
any $X\in\Ii^+_\sigma(\omega) \A(\rho')$
has an analytic extension on $T_a\DD(\rho)$. 
The linear map $\TT_a(X)= T_a^{-1}\,X\circ T_a$ from $\Ii^+_\sigma(\omega)\A(\rho')$
to $\A'(\rho)$ is compact with norm $\|\TT_a\|\leq 6\pi
a/(\rho'-\kappa\rho) + 3a$. 
\end{proposition}

The proof of this proposition is contained in Section \ref{section
sets I+ and Ik}. 
As it will be seen, 
it is the existence of $\kappa<1$ such that 
 $I^+_\sigma(\omega) \subset \{k\in\Zz^2\colon \|T^*_ak\|\leq \kappa \|k\|\}$
that guarantees the above analyticity improvement.

The following theorem (proved in \cite{jld} using a ``homotopy method'') states the existence of a nonlinear change of coordinates $U$ isotopic to the
identity, eliminating all the far from resonance terms $I^-_\sigma(\omega')=I^-_{\sigma/\alpha'}(\omega'/\alpha')$ of
a vector field $X$ in a neighbourhood of a non-zero constant vector field.

\begin{theorem}[\cite{jld}]\label{main theorem1}
Let $0<\rho'<\rho$, $\psi\in \Cc^d\setminus\{0\}$ and
$0<\sigma<\|\psi\|$.
If $X$ is in the open ball $\hat B\subset\A'(\rho)$ centred at
$\psi$ with radius 
$$
\hat\varepsilon=
\frac{\sqrt6-2}{12}\sigma
\min\left\{
\frac{\rho-\rho'}{4\pi},
\frac{3-\sqrt6}{6} \frac\sigma{\|\psi\|}
\right\},
$$ 
there is a diffeomorphism $U\colon\DD(\rho')\to\DD(\rho)$
isotopic to the identity, satisfying 
$$
\begin{array}{c}
[\Ii_\sigma^-(\psi)](DU)^{-1}\,X\circ U=0,
\quad\text{ and }\quad
U=\id 
\text{ if }[\Ii_\sigma^-(\psi)]X=0.
\end{array}
$$ 
The map $\UU_\psi\colon \hat B \to [\Ii_\sigma^+(\psi)]\A(\rho')$ given by $X\mapsto (DU)^{-1}\,X\circ U$
is analytic, and the derivative at $\psi$
is $\Ii_\sigma^+(\psi)$.
Moreover, 
$$
\|\UU_\psi(X)-\psi\|_{\rho'}\leq
2\left(
1+\max\left\{
\frac23(3-\sqrt6),6 (\sqrt6+2)\frac{\|\psi\|}\sigma
\right\}
\right)
\|X-\psi\|'_\rho.
$$ 
\end{theorem}

The definition below uses a different
order of the linear and nonlinear maps, compared with the approach of
\cite{jld}. 
There is no fundamental difference, only a technical adjustment to
allow us to evaluate the size of the iterates right after performing
$\UU$.

\begin{definition}
A map $\RR_\omega$ acting on vector fields $X$ is
called a {\it one-step renormalisation} operator\index{renormalisation operator} if it is of the form
$$
\RR_\omega(X)=\frac{1}{\hat\omega'\cdot\Ee(X')}\, X', \quad \text{with}\quad
X'=\UU_{\omega'/\alpha'}\circ\TT_a(X).
$$
The nonlinear transformation 
$\UU_{\omega'/\alpha'}$ that eliminates the far from resonance terms $I^-_\sigma(\omega')$, is defined according to Theorem \ref{main theorem1} with respect
to $\omega'/\alpha'$. 
The linear map $\TT_a$ is given by Proposition \ref{proposition L}.
\end{definition}

Since $\RR_\omega(\omega)=\omega'$,
fixed points of $\RR_\omega$ are possible for constant vectors $\omega$
such that $\omega'=\omega$. Those correspond to the ones with
slope $\alpha$ having a constant continued fraction expansion, i.e. of the form
$[a, a,\dots]$, $a\geq1$. 
Other periodic points of $\RR_\omega$ can be obtained for $\alpha$ with a periodic continued
fraction expansion, which, for period $p$, is of the form $[\overline{a_0, \dots, a_p}]$,
$a_0,\dots,a_p\geq1$.
Eventually periodic points appear for slopes of vectors with eventually
periodic continued fraction expansion (see \cite{jld}).

\begin{proposition}\label{proposition existence of one-step renormalisation}
Given $\rho'>0$ and $\sigma>0$ sufficiently small, we can find $C>0$
such that the one-step
renormalisation operator $\RR_\omega$ is a well-defined analytic
map from an open ball $B$ of $\Ii^+_\sigma(\omega)\A(\rho')$ centred
at $\omega$ with radius $\zeta=C\sigma^2/(a\|\omega'\|)$, to
$\Ii^+_\sigma(\omega')\A(\rho')$.
In addition, $\|\RR_\omega(\omega+f)-\omega'\|_{\rho'}\leq
\|f\|_{\rho'}/\zeta$ and
$\|\RR_\omega(\omega+f)-\omega'-D\RR_\omega(\omega)\,f\|_{\rho'} \leq
[\zeta(\zeta-\|f\|_{\rho'})]^{-1}\|f\|^2_{\rho'}$.
\end{proposition}

\proof
Suppose $\rho>\rho'$ and $\kappa<1$ such that $\kappa\rho<\rho'$.
We will always make these choices such that $\rho'-\kappa\rho$ is
bounded away from zero.
Then, by Proposition \ref{proposition L}, $\TT_a(B) \subset B'$, with $B'=\{X\in \A'(\rho) \colon
\|X-\omega'/\alpha'\|'_\rho< 
CC'\sigma^2/\|\omega'\|\}$, where $C'$ is a
positive constant.
Then, $C$ has to be chosen sufficiently small so that
the transformation $\UU_{\omega'/\alpha'}$ is valid in a neighbourhood of
$\omega'/\alpha'$ containing $B'$,
according to Theorem \ref {main theorem1}.
The same theorem states that
$\UU_{\omega'/\alpha'}(B')\subset
B''$ where
$B''=\{X\in\Ii^+_{\sigma/\alpha'}(\omega'/\alpha')\A(\rho')\colon\|X-\omega'/\alpha'\|_{\rho'}<CC'C''\sigma/\alpha'\}$,
where $C''>0$ is an independent constant.

Let $F$ be the complex-valued continuous functional
$F(X')=\hat\omega'\cdot\Ee(X')$,
with the domain being a neighbourhood
$B'''\subset\A(\rho')$ of $\omega'/\alpha'$ such that $B''\subset
B'''$ and
$F(B''')\subset\{z\in\Cc\colon|\alpha'z-1|<1/2\}$. 
Writing $f'=X'-\omega'/\alpha'$, it is enough to have
$\|f'\|_{\rho'}<\ell(4\sqrt2-5)/2\leq1/(2\alpha'\|\hat\omega'\|)$ with
$\omega'=\ell(1,\alpha')$,
which is satisfied for a sufficiently small $C$.
Hence
$[\hat\omega'\cdot\Ee(X')]^{-1}$ is bounded and analytic in
$B''$ and, for $X'\in B''$,
$\|X'/(\hat\omega'\cdot\Ee\,X')-\omega'\|_{\rho'}\leq
CC'C''C'''\sigma$, for some independent scalar $C'''>0$.
This, Theorem \ref{main theorem1} and Proposition
\ref{proposition L} prove the first part of the claim.
The first estimate also follows in the same way.

To show the estimate on the second order remaining of the Taylor
expansion of $\RR_\omega$ around $\omega$, consider $g\colon z\mapsto
\RR_\omega(\omega+zf)-\omega'$, $f\in B$.
This is an analytic map on an open ball in $\Cc$ with radius
$r=\zeta/\|f\|_{\rho'}>1$. So, Cauchy's formula gives
\begin{eqnarray*}
\|g(1)-g(0)-g'(0)\|_{\rho'}
&\leq&
\frac1{2\pi}\oint_{|z|=r}\frac{\|g(z)\|_{\rho'}}{|z^2(z-1)|}dz\\
&\leq&
\frac1{ r(r-1)}\sup_{\|\xi\|_{\rho'}=\zeta}\|\RR_\omega(\omega+\xi)-\omega'\|_{\rho'}.
\end{eqnarray*}
\cqd

%%%%%%%%%%%%%%%%%%%%%%%%%%%%%%%%%%%%%%%%%%%%%%%%%%%%%%%%%%%%%%%%%%%%%%%%%%%%

\section {Diophantine Numbers and $GL(2,\Zz)$}\label{Diophantine Numbers}

Let us review the definition of diophantine numbers and some of
its properties, and introduce the set of vectors of the plane with
diophantine slope.

\begin{definition}
For an irrational $\alpha$ and $\beta\geq0$, we say
that $\alpha$ is a {\it diophantine number of order $\beta$}\index{diophantine!number} if there is a
constant $C>0$ such that 
$$
\left|\alpha-\frac pq\right| > \frac {C}{q^{2+\beta}}
$$
for any $p/q\in\Qq$. We denote the set of all the diophantine
numbers of order $\beta$ by $DC(\beta)$.
All non-zero vectors $\omega$ in $\Rr^2$ with slope $\alpha\in
DC(\beta)$ form the set denoted also by
$DC(\beta)$, which is $GL(2,\Zz)$-invariant.
This is equivalent to say that we can find $C>0$ such that
$|\omega\cdot k| >C\|k\|^{-(1+\beta)}$, $k\in\Zz^2\setminus\{0\}$.
\end{definition}

If $\alpha\in DC(\beta)$ and using (\ref{bounds on (alpha-p/q) and beta_n}), there exists $K>1$ such that
$$
\frac1{q_nq_{n+1}}>\abs{\alpha-\frac{p_n}{q_n}}>\frac
{K^{-1}}{q_n^{2+\beta}}, \quad n\geq0,
$$
where $p_n$, $q_n$ are the convergents of $\alpha$. 
This yields, together with (\ref{recurrence
formulae q_n, p_n}), (\ref{upper and lower bounds on beta_n}) and (\ref{upper and lower bounds on tilde A(n)}), that
the diophantine condition can be written in the following ways:
\begin{equation}\label{diophantine conditions}
\begin{array}{rl}
q_{n+1} & <Kq_n^{1+\beta},
\\ 
a_{n+1} & <Kq_n^\beta,
\\
\beta_{n+1}^{-1} & < 2K \beta_n^{-(1+\beta)},
\\
\tilde A_{n+1} & < 2K \tilde A_n^{1+\beta}, \quad n\geq0.
\end{array}
\end{equation}

Finally, for $\beta>0$, one obtains lower and upper bounds
for the convergents
$F_n \leq q_n <  K^{[(1+\beta)^n-1]/\beta}$,
and for the coefficients $a_{n+1}<K^{(1+\beta)^n}$, where $F_n$ is
the Fibonacci sequence \index{Fibonacci sequence} given by the recurrence formula
$F_{n+1}=F_n+F_{n-1}$ (i.e. with $q_n=F_n$ and $a_n=1$).
The sequence $F_n$ has an exponentially growing solution of the
form $\gamma^n$. 
All the {\it constant type numbers}\index{constant type numbers}, i.e. the ones inside $DC(0)$,
have an exponential bound on the sequence $q_n$ similar to $F_n$, and $a_n<K$.

The set of the vector fields $X$ close enough to $\Ee(X)$ with winding
ratio\index{winding ratio}
$w\in DC(\beta)$ is $\RR_\omega$-invariant for any $\omega\in DC(\beta)$. 
This is easily seen by recalling that 
if $X'$ is the new vector field after a coordinate transformation
isotopic to $T\in GL(2,\Zz)$, then the winding ratio of $X'$ is $Tw/\|Tw\|$.
In particular, $DC(\beta)$ is $T_a$-invariant, and $T_a$ is the linear map
involved in $\RR_\omega$, modulo a constant rescaling.

The constant type numbers $DC(0)$ have zero Lebesgue measure, but
the {\it Roth numbers} (see e.g. \cite{Cassels,Herman}): $\cap_{\beta>0} DC(\beta)$, have
measure one, which implies the full measure of any set
$DC(\beta)$, $\beta>0$.
Note that if $\beta_1<\beta_2$ then
$DC(\beta_1)\subset DC(\beta_2)$.

%%%%%%%%%%%%%%%%%%%%%%%%%%%%%%%%%%%%%%%%%%%%%%%%%%%%%%%%%%%%%%%%%%%%%%%%%%%%

\section {The Limit Set of the Renormalisation}\label{The Limit Set of the Renormalisation}

Let $\omega_0\in\Rr^2\setminus\{0\}$ with non-zero slope
$\alpha_0=[a_0,a_1,\dots]$.
For an initial vector field $X_0$, define the sequence of its
images under renormalisation as $X_{n+1}=\RR_{\omega_n}(X_n)$,
whenever possible for $n\geq0$.
We denote also $\omega_{n+1}=\RR_{\omega_n}(\omega_n)$.
Notice that, if $w(X_n)=w(\omega_n)$, then 
$w(X_{n+1})=w(\omega_{n+1})$.

The dynamics of $\omega_n$ with
slope equal to $\alpha_n=[a_n,a_{n+1},\dots]$, is given by the shift
of coefficients, related to the Gauss map.
If $\alpha_0$ is a rational number, then its continued fraction
expansion is finite.
Therefore, it is not possible to perform the
iteration of this procedure for an infinite number of times.

We want to include in the domain of this renormalisation iteration, vectors
$\omega_0$ with any slope and vector fields $X_0$ with also far
from resonance terms.
So, initially, we allow $\RR_{\omega_0}$ to be an ``adjustment'' procedure,
using the change of basis $V$ or $S$ and $\UU$.
This unique transient step does not influence the qualitative result
that follows.

\begin{theorem}\label{theorem generic frequency d=2}
Let $\rho'>0$.
For any vector $\omega_0\in \cup_{\beta<1} DC(\beta)$,
there exists an open neighbourhood $B$ of $\omega_0$ in $\A(\rho')$
such that: if $X_0\in B$ and $w(X_0)=w(\omega_0)$, then 
$$
\| X_n - \omega_n \|_{\rho'}  <  K \theta ^ {n}, \quad n\geq0,
$$
where $K>0$ and $0<\theta<1$.
\end{theorem}

The proof of this theorem is included in Section
\ref{proof of theorem generic frequency d=2}.
For ``irrational winding ratios'' not in the condition of Theorem \ref{theorem generic frequency d=2}, the domain of the
renormalisation shrinks faster than
the rate of convergence of the iteration.
Both of these behaviours are determined by the size of the coefficients $a_n$.
If their growth is controlled as for the diophantine case, the renormalisation reduces the perturbation enough to keep
the iterates inside the domain.

A corollary of the above theorem relates to the main result 
in \cite{jld} for the two dimensional case.
In fact, quadratic irrationals have an eventually periodic continued fraction expansion
and are inside $DC(0)$.
It is then enough to choose a renormalisation operator $\RR$
associated to such $\omega$, with
the block of matrices $T$ corresponding to the periodic string.
In this way, $\omega$ is an eigenvector of $T$ and a fixed point of $\RR$.
Theorem \ref{theorem generic frequency d=2} claims the existence of an
invariant set contracting towards $\omega$,
that includes all vector fields in some neighbourhood with
same winding ratio.
To show that that is a submanifold, and there is a
one-parameter family corresponding to the unstable direction,
we refer to the details of the proof (namely the spectral properties
of the derivative of $\RR$ at $\omega$).
The fixed point is therefore hyperbolic.

Another application of Theorem \ref{theorem generic frequency
d=2} is a result for the Poincar\'e map on the circle, in the same
spirit as the one proved
originally by Arnol'd \cite{Arnold2} in the analytic case:
any real analytic circle map with diophantine rotation number
$\alpha$, inside some neighbourhood of the rotation
$R_\alpha\colon x\mapsto x+\alpha \mod 1$,
is analytically conjugate to $R_\alpha$, for a norm
induced by the one used here for vector fields
(cf. Moser \cite{Moser2} for the differentiable case, and
Herman \cite{Herman} and Yoccoz \cite{Yoccoz} for more general and global results).
We have not determined the optimal condition obtained by Yoccoz
\cite{Yoccoz2} in the context of the linearisation of germs of
analytic diffeomorphisms and the local
theorem on analytic conjugacy of circle
diffeomorphisms.
That condition states that it is necessary and sufficient that the set
of winding ratios for which the claim in Theorem \ref{theorem generic frequency
d=2} is expected to be valid strictly contains
$\cup_{\beta\geq0} DC(\beta)$, with slope verifying the Brjuno
condition \index{Brjuno condition}:
$\sum_{n\geq0}\log (q_{n+1}) /q_n<+\infty$.

%%%%%%%%%%%%%%%%%%%%%%%%%%%%%%%%%%%%%%%%%%%%%%%%%%%%%%%%%%%%%%%%%%%%%%%%%%%%

\section{The Sets $I^+_\sigma(\omega)$ and $I^\kappa_a$} \label{section sets
I+ and Ik}

Here it is enough to consider
$\omega=(\omega_1,\omega_2)$, $\omega_1,\omega_2>0$, with slope
$\alpha=\omega_2/\omega_1 > 1$, and $a=[\alpha]$.

\begin {definition}\label{definition Ik}
For $\kappa>0$ define the set
$$
I^\kappa_a=\{k\in\Zz^2\colon \|T_a^*k\|\leq\kappa\|k\|\},
$$
and the projection $\Ii^\kappa_a$ on $\A(r)$, $r>0$, given by
$$
\Ii^\kappa_a X(\theta)=\sum\limits_{k\in I^\kappa_a}X_ke^{2\pi i k
\theta}, \quad X\in\A(r).
$$
\end {definition}

Let $0<\sigma<\omega_i$, $i=1,2$, and $1/2<\kappa<1$.
From the respective definitions one can rewrite the form of the cones $I^+_\sigma(\omega)$ and
$I^\kappa_a$, closed subsets of $\Zz^2$.
In fact,
$I^+_\sigma(\omega)$ is
bounded by the lines $k_2=mk_1$ and $k_2=lk_1$, and
$I^\kappa_a$ by $k_2=sk_1$ and $k_2=rk_1$,
where 
\begin{equation}\label{bounds on slopes}
\begin{array}{ccccc}
m=-\frac{\omega_1-\sigma}{\omega_2+\sigma}, &
l=-\frac{\omega_1+\sigma}{\omega_2-\sigma} & \text{ and } & 
s=-\frac{1-\kappa}{a-1+\kappa}, &
r=-\frac{1+\kappa}{a+1-\kappa}.
\end{array}
\end{equation}
Note that $l<m<0$ and $r<s<0$ if $\kappa>(a+1)^{-1}$.

\begin{lemma}\label{lemma R(I+) contained Ik}
If 
$\kappa \geq 1-
\|\omega\|^{-1}
\min \{
a(\omega_1-\sigma),
2(\omega_2-\sigma)-a(\omega_1+\sigma)
\}$,
then $I^+_\sigma(\omega)\subset I^\kappa_a$.
Moreover, considering $\omega=(1,\alpha)$ with $\alpha>1$ and
$\sigma<1/3$, it is sufficient to have
$\kappa\geq 1-\eta$ where $\eta=(1-3\sigma)/3>0$.
\end{lemma}

\proof
It is enough to check the values of $\kappa$ for which we have $r\leq
l\leq m \leq s$, as given in (\ref{bounds on slopes}).
The result follows from a simple calculation, considering the conditions on $\sigma$
imposed before.
\cqd

%%%%%%%%%%%%%%%%%%%%%%%%%%%%%%%%%%%%%%

\subsection{Analyticity Improvement in
$I^\kappa_a$}\label{Analyticity Improvement in Ik}

In order to prove that Proposition \ref{proposition L} holds, we need
to show that $X\circ T_a$ is analytic in
$\DD(\rho)$ and has bounded derivative.
From the result of Lemma \ref{lemma R(I+) contained Ik}
it is enough to notice that, for $X(\theta)=\sum\limits_{I^\kappa_a} f_k
e^{2\pi i k\cdot\theta}$,
$$
\|X\circ T_a\|_\rho  \leq\sum\limits_{k\in I^\kappa_a}\|f_k\|e^{\rho\|T_a^*k\|} 
 \leq \sum\limits_{k\in I^\kappa_a} \|f_k\| e^{(\rho\kappa-\rho')\|k\|}e^{\rho'\|k\|} 
\leq \|X\|_{\rho'},
$$
and
$$
\|D(X\circ T_a)\|_\rho  \leq 2\pi \sum\limits_{k\in I^\kappa_a}\|T_a^*k\|
e^{-\delta\|k\|} \|
f_k\|e^{(\rho\kappa+\delta-\rho')\|k\|}e^{\rho'\|k\|} 
\leq \frac{2\pi}{\delta}\kappa\|X\|_{\rho'},
$$
by choosing $0<\delta<\rho'-\kappa\rho$, and making use of the relation
$\sup_{t\geq0}te^{-\xi t}\leq 1/\xi$ for $\xi>0$.
These bounds imply that $\| \TT_a(X) \|'_\rho \leq
(1+{2\pi\kappa}/\delta)\|T_a^{-1}\|\,\|X\|_{\rho'}$. 
We choose e.g. $\delta=\kappa(\rho'-\kappa\rho)$.

The above for $\DD(\rho)$ is true also for $\DD(r)$, $r>\rho$ such that $\rho'>r\kappa$.
Therefore, $\TT_a=\II\circ\JJ$, where
$\JJ\colon\Ii^\kappa_a\A(\rho')\to\A'(r)$ is bounded as $\TT_a$ above, and the inclusion map
$\II\colon\A'(r)\to\A'(\rho)$,
$\II(X)=X|_{\DD(\rho)}$ is
compact.

%%%%%%%%%%%%%%%%%%%%%%%%%%%%%%%%%%%%%%%%%%%%%%%%%%%%%%%%%%%%%%%%%%%%

\section{Convergence Result}
\label{proof of theorem generic frequency d=2}

Here we include the proof of Theorem \ref{theorem generic frequency
d=2}, the convergence of the iterates $X_n$ towards the projected
space $\Ee\,\A(\rho')=\Cc^2$ of constant vector fields.
More precisely, we show that the orbit of $X_0$ approximates
exponentially the orbit of $\omega_0$.
For all steps of the iteration, we choose the same $\kappa$ and $\rho$
from Proposition \ref{proposition L}.

Considering the $n$-th step, 
the slope of $\omega_n$ is given by
$\alpha_n>1$, and we denote its integer part by $a_n=[\alpha_n]$.
While $\alpha_n>1$ for $n\geq1$, we can have 
$\alpha_0>0$, but this can be regarded as a transient step (using $S$) which does
not carry any problems into the following analysis.
Also, if $\alpha_0$ is negative, 
the change of basis $V$
can be used initially, returning to the case of positive slope.
Therefore, we will assume here $\alpha_n>1$, $n\geq0$.
Again, the same can be said about starting within a space of vector
fields with more than resonant terms. 
After an initial transformation, the iterates will always be inside a
space restricted to vector fields with no far from resonance terms.

Let $\sigma>0$ as given in Proposition \ref{proposition L}.
To simplify the notation
replace $\RR_{\omega_n}$ simply by $\RR_n$ and use
the projection $\Pp_{n}$ of $\A(\rho')$ over $\omega_n$:
$$
\Pp_{n}f= (\hat\omega_{n}\cdot f)\,\omega_{n}, \quad
f\in\A(\rho'),
$$ 
where $\hat v=v/(v\cdot v)$, for every $v\in\Rr^2\setminus\{0\}$.
According to Proposition \ref{proposition existence of one-step
renormalisation}, the domain of each $\RR_n$ is the ball
$B_n\subset\A(\rho')$ around $\omega_n$ with radius $\zeta_n$.
The derivative at $\omega_n$ is the linear operator
given by:
$$
D\RR_n(\omega_n)=
(\Ii - \Pp_{n+1}\Ee)\,\LL_n,
$$
with
$\LL_n = \alpha_{n+1}\Ii^+_{\sigma}(\omega_{n+1}) \circ\TT_{a_n}$.
This is a compact operator since $\Ii^+$ is bounded and $\TT_{a_n}$ is compact.

The following proposition shows that there is a ``super''
exponentially shrinking of the non-constant terms, the image of $(\Ii-\Ee)$,
corresponding to the zero eigenvalue of the stable subspace.
Recall that $\tilde A_n=\prod_{i=0}^n \alpha_i$ and $\alpha_0q_n <
\tilde A_n < 2\alpha_0 q_n$, with $q_n$ growing at least
exponentially with $n$.
Also, let $\tilde A_{-1}=1$.

\begin{proposition}\label{convergence on stable mfld}
Suppose $\omega_0=\ell(1,\alpha_0)\in DC(\beta)$, $\beta\geq0$,
$\ell\not=0$.
If $\sigma<\frac\ell 2(1+4\gamma^3)^{-1}$, then
there exists $c_1,c_2>0$ such that
$$
\|\LL_{n}\circ \dots\circ\LL_j \,(\Ii-\Ee)\| \leq
c_1 e^{ -c_2 \Lambda_{j,n} } \|\LL_{n}\circ \dots\circ\LL_{j+1} \,(\Ii-\Ee)\|,
$$ 
for $j=0,\dots,n$ with $n>0$, and where
$$
\Lambda_{j,n}^{2+\beta}=\frac{
\tilde A_{n+1}\tilde A_n
}{\sigma \tilde A_{j-1}^{1+\beta} \tilde A_{j-1}}.
$$
\end{proposition}

\proof
Consider a vector $\omega$ with slope $\alpha$, $a=[\alpha]$, the
matrix $T_a$ and $\lambda=(a+\sqrt{a^2+4})/2>a$. 
Let $\Omega$ be an orthogonal vector to $\omega$,
e.g. $\Omega=(1,-1/\alpha)$.
So, 
\begin{equation}\label{iteration of w and Omega}
T_a^{-1}\omega=\frac1{\alpha'} \omega' \quad\text{ and } \quad
T_a\Omega=-\frac1{\alpha}\Omega',
\end{equation}
where 
$\Omega'=(1,-1/{\alpha'})$.
The projections of $k\in I^+_{\sigma}(\omega)$ over $\omega$ and
$\Omega$ are $\Pp_\omega \,k=(\omega\cdot k)\,\hat\omega$ and
$(\Ii-\Pp_\omega)\,k=(\hat\Omega\cdot k)\,\Omega$, respectively.

Introducing the notation $I^+_n=I^+_{\sigma}(\omega_n)\setminus\{0\}$,
$n\geq0$, and the corresponding projection $\Ii^+_n$, define a subset
$V^+_{j,n}$ 
of $I^+_j$ as 
$$
V_{j,n}^+=\{k\in I^+_j \colon
T^*_{i}\cdots T^*_j \, k \in I^+_{i+1}, j\leq i\leq n \},
$$
for $j=0,\dots,n$ when $n>0$.
This set includes all the Fourier modes in $I^+_j$ that will be mapped
into $I^+_{n+1}$ (the only relevant ones since $\LL_n$ includes the
projection $\Ii^+_{n+1}$) by the sequence of matrices $T^*_{n}\cdots T^*_j$. 
Therefore, for every $k$ in $V^+_{j,n}$, we can write
$|\omega_{n+1}\cdot T^*_{n}\dots T^*_j \, k| \leq \sigma \|T^*_{n} \dots
T^*_j \, k\|$.
Notice that $|\omega_{n+1}\cdot T^*_{n}\dots T^*_j \, k| = |T_j\dots
T_{n}\, \omega_{n+1}\cdot k|$ and, from (\ref{iteration of w and Omega}),
$$
T_j\dots T_{n} \, \omega_{n+1}= \left(\prod_{i=j+1}^{n+1}\alpha_i\right) \,
\omega_j.
$$
It follows that 
\begin{equation}\label{inequality for I_j}
|\omega_j \cdot k| \prod_{i=j+1}^{n+1}\alpha_i  \,
 \leq \sigma \|T^*_{n} \dots
T^*_j \, k\|, \quad k\in V^+_{j,n}.
\end{equation}
A formula for the calculation of $T^*_{n} \dots T^*_j \, k$, with $k\in
V^+_{j,n}$, can be
obtained by induction:
\begin{equation}\label{formula for T*n ... T*j k}
\begin{array}{rl}
T^*_{n} \dots T^*_j \, k = & (\omega_j\cdot k)\left[\hat\omega_{n+1}
\prod\limits_{i=j+1}^{n+1}\alpha_i  
\right.
\\
&
+
\left.
\Omega_{n+1} 
\sum\limits_{m=j}^{n}(\hat\Omega_{m+1}\cdot T^*_m\,\hat\omega_m)
{\left(\prod\limits_{i=j+1}^{m} \alpha_i\right)}
{\prod\limits_{i=m+1}^{n}(-\alpha_{i})^{-1}}\right] \\
& + 
\Omega_{n+1}\,(\hat\Omega_j\cdot k)\prod\limits_{i=j}^{n}(-\alpha_i)^{-1},
\end{array}
\end{equation}
making use of the orthogonality between $\omega_i$ and $\Omega_i$ for any $i$,
i.e. $\omega_i\cdot\Omega_i=0$.
That yields
$T^*_i\,\hat\omega_i= (\hat\Omega_{i+1}\cdot
T^*_i\,\hat\omega_i)\,\Omega_{i+1} + (\omega_{i+1}\cdot T^*_i\,\hat\omega_i)\,
\hat\omega_{i+1}$, with
$\omega_{i+1}\cdot T^*_i\,\hat\omega_i=T_i\,\omega_{i+1}\cdot
\hat\omega_i=\alpha_{i+1}$. 
%Using (\ref{bound on x by K k}), 
There is $K>1$ such that
the norm of $T^*_{n} \dots T^*_j \, k$ as given in (\ref{formula for
T*n ... T*j k}) is bounded by
$$
\begin{array}{rl}
\|T^*_{n} \dots T^*_j \, k \| \leq
&
|\omega_j\cdot k| \left[ \|\hat\omega_{n+1}\| \,
\prod\limits_{i=j+1}^{n+1}\alpha_i \, 
\right.
\\
&
\left.
+
2 \sum\limits_{m=j}^{n}
|\hat\Omega_{m+1}\cdot T^*_m\,\hat\omega_m|
\left(\prod\limits_{i=j+1}^{m} \alpha_i\right)
\prod\limits_{i=m+1}^{n}\alpha_{i}^{-1}
\right]
\\
&
+
2K\left(\prod\limits_{i=j}^{n}\alpha_i^{-1}\right) \, \|k\|, \quad k\in V^+_{j,n}.
\end{array}
$$
We now need an upper estimate of the summation in the above inequality.
If $\alpha$ is a quadratic irrational ($\omega$ and $\Omega$ are
orthogonal eigenvectors of $T$), then $\hat\Omega_{m+1}\cdot
T^*_m\,\hat\omega_m=0$ for every $m$, 
which is in agreement with \cite{jld}.
On the other hand, for other values of $\alpha$, notice that
$|\hat\Omega_{m+1}\cdot
T^*_m\,\hat\omega_m|<\frac1\ell(\alpha_{m+1}^{-1}+\alpha_m^{-1})<\frac2\ell$.
In addition,
$$
\begin{array}{rl}
\sum\limits_{m=j}^{n}\left(\prod\limits_{i=j+1}^{m}\alpha_i\right)\prod\limits_{i=m+1}^{n}\alpha_i^{-1}
&
\leq
\alpha_{n+1}^{-1}\prod\limits_{i=j+1}^{n+1}\alpha_i \sum\limits_{m=j}^{n}\tilde
A_m\tilde A_n^{-1}
\\
&
\leq
\gamma^3\alpha_{n+1}^{-1} \prod\limits_{i=j+1}^{n+1}\alpha_i,
\end{array}
$$
because, from the second inequality in (\ref{bound on tilde A(n)}), we
have that 
$$
\sum_{m=j}^{n}\tilde A_m\tilde A_n^{-1}\leq\gamma^3.
$$
The hypothesis and the formula
$\omega_n=\ell(1,\alpha_n)$, where $\omega_0=\ell(1,\alpha_0)$, imply that
$\Gamma_n\sigma<1/2$ with
$\Gamma_n=\|\hat\omega_n\|+4\gamma^3/(\alpha_{n}\ell)$.
Then, using (\ref{inequality for
I_j}) it follows that
\begin{equation}\label{relation for w_n.k a}
|\omega_j\cdot k| \leq \frac{2K\sigma\|k\|}{1-\Gamma_{n+1}\sigma}
\prod_{i=j}^{n}(\alpha_i\alpha_{i+1})^{-1} 
\end{equation}

Given $j\geq0$, as $\omega_0\in DC(\beta)$,
there is a constant $C_0>0$ such that
\begin{equation}\label{result of diophantine cdn a}
|\omega_j\cdot k|= |\omega_0\cdot
T_{0}^{*-1} \dots T_{j-1}^{*-1} \, k| 
\prod_{i=1}^{j}\alpha_i \,
\geq \frac{C_0 \prod_{i=1}^{j}\alpha_i}{\|T_{0}^{*-1} \dots T_{j-1}^{*-1} \, k\|^{1+\beta}}.
\end{equation}
Consider the supremum norm on $\Zz^2$,
$\|k\|_\infty=\max_{i=1,2}|k_i|$.
It is equivalent to $\|\cdot\|$ since
$\frac12\|\cdot\| \leq \|\cdot\|_\infty \leq \|\cdot\|$.
Then, $\| T_{0}^{*-1} \dots T_{j-1}^{*-1} \, k\|_\infty 
=\|P^{-1}_{j-1}k\|_\infty=\max\{p_{j-1}+p_{j-2},q_{j-1}+q_{j-2}\}
\leq
b\tilde A_{j-1} \, \|k\|_\infty$, for some $b>0$, implies that
\begin{equation}\label{bound on T-1...k a}
\| T_{0}^{*-1} \dots T_{j-1}^{*-1} \, k\| \leq 2b\tilde A_{j-1} \, \|k\|.
\end{equation}

Now, it follows from inequality (\ref{relation for w_n.k a}) together with 
(\ref{result of diophantine cdn a}) and
(\ref{bound on T-1...k a}) that, for all $k\in V^+_{j,n}$,
$$
\|k\|^{2+\beta} \geq
\frac{C}{\sigma 
\tilde A_{j-1}^{1+\beta}} \left(\prod_{i=1}^{n+1}\alpha_i\right) \,
\prod_{i=j}^{n}\alpha_i 
=
\frac{C \tilde A_n\tilde A_{n+1}}{\sigma \tilde A_{j-1}^{1+\beta}\tilde A_{j-1}},
$$
with $C=C_0/(2^{3+\beta}\alpha_0b^{1+\beta}K)$.

Let the 
operator $\Vv^+_{j,n}\colon\A(r)\to\A(\rho')$, $r>\rho'$, be a
projection over the indices in $V^+_{j,n}$ together
with an analytic inclusion.
That is,
$$
\|\Vv^+_{j,n}f\|_{\rho'}  = \sum\limits_{k\in V^+_{j,n}}\|f_k\|
e^{r\|k\|}e^{-(r-\rho')\,\|k\|}  \leq \|\Vv^+_{j,n}\|\, \|f\|_r,
$$ 
where
$\|\Vv^+_{j,n}\| \leq \exp\left[-C^{1/(2+\beta)}(r-\rho') \Lambda_{j,n}\right]$.
It is possible to include this transformation in the calculation
of the sequence
$$
\LL_n\circ\dots\circ\LL_{j}\,(\Ii-\Ee)=\LL_n\circ\dots\circ\LL_{j+1}
\,(\Ii-\Ee)\,\Vv^+_{j,n}\bar\LL_{j},
$$
$0\leq j\leq n$,
where $\bar\LL_j\colon\A(\rho')\to\A(r)$ with $\rho'>\kappa r$, is $\LL_j$ followed
by an analytic extension.
The norm of $\bar\LL_j$ can be estimated above by $\alpha_{j+1}\alpha_j$ up to
the product with a constant.
Therefore,
$$
\|\LL_n\circ\dots\circ\LL_{j}\,(\Ii-\Ee)\|  
\leq
c_1 \exp (-c_2 \Lambda_{j,n})
\|\LL_n\circ\dots\circ\LL_{j+1}\,(\Ii-\Ee)\|
$$
with $c_1, c_2$ some positive constants independent either of $j$ or $n$.
\cqd

To complete the description of the eigenspaces of $D\RR_n(\omega_n)$,
we look at the constant terms, the image of $\Ee$.
That is, we have to solve
$G_n(f)= \lambda f$, with
$$
G_n(f) =
D\RR_n(\omega_n)\Ee\,(f)=(\Ii-\Pp_{n+1})
\alpha_{n+1}T^{-1}_{n}\Ee\,(f).
$$
This operator can be re-written in the following way:
$$
G_n(f)=\frac{\alpha_{n+1}}{1+\{\alpha_n\}^2}
\left[
\begin{matrix}
-\alpha_n & 1 \\
\{\alpha_n\}\alpha_n & -\{\alpha_n\} \\
\end{matrix}
\right]
f_0.
$$
As it has zero determinant its eigenvalues are 0 and $\nu_n=\tr(G_n)$,
with $|\nu_n|>1$ for any $\alpha_n>0$.
The respective eigenvectors are $\omega_n$ and
$\Omega_{n+1}$, which determine the stable and unstable
subspaces corresponding to the constant term.

The constant part of each iterate $X_n=\omega_n + f_n$ can be written as
$\Ee(X_n)=\omega_n + \Pp_{n}\circ\Ee\,(f_{n}) + (\Ii-\Pp_{n})\circ\Ee\,(f_{n})$.
Notice that
\begin{equation}\label{P E f_n+1}
\begin{array}{rl}
\Pp_{n+1}\circ\Ee\,(f_{n+1}) & 
= \Pp_{n+1}\circ G_n \, (f_n) + \OO_n(\|f_n\|_{\rho'}^2) \\
&
= \OO_n(\|f_n\|_{\rho'}^2),
\end{array}
\end{equation}
as we know from above that the image of $G_n$ is in
$(\Ii-\Pp_{n+1})\,\Cc^2$.
That is, its size can be estimated by the square of the norm of $f_n$.
The other term, $(\Ii-\Pp_{n})\circ\Ee\,(f_{n})$, can be controlled in a different way, as
it is possible to exclude some vector fields in $B_n$ that do not
have winding ratio $\omega_n/\|\omega_n\|$.

\begin{lemma}\label{lemma of the winding ratio cone}
Suppose that $X\in B_n$ has winding ratio\index{winding ratio} $\omega_n/\|\omega_n\|$.
Then $X$ belongs to the subset
$$
C_n=\left\{ X\in B_n\colon 
\|(\Ii-\Pp_n)\circ\Ee\,(X) \|
\leq
\| (\Ii-\Ee) X\|_{\rho'}
 \right\}.
$$
\end{lemma}

\proof
A subset of vector fields $D_n\subset B_n$ that do not cross the line
spanned by $\omega_n$ can be of the form:
$$
D_n=\left\{ X\in B_n\colon \| X(\theta)-\Ee\,X \| <
\|(\Ii-\Pp_n)\circ\Ee\,(X)\| , \theta\in\DD(\rho')\right\}.
$$
That is, the
non-constant part of $X\in D_n$ at each point $\theta$ is less
than the distance (given by the norm $\|\cdot\|$) between
$\Ee\,(X)$ and its projection by $\Pp_n$ over the subspace spanned by
$\omega_n$.
The slopes of all the vectors $X(\theta)$ are bigger than $\alpha_n$ or
always less than $\alpha_n$, never crossing that value (as for their
respective winding ratios).
Therefore,
since we have that 
$\| X(\theta)-\Ee\,X \| \leq \| (\Ii-\Ee) X\|_{\rho'}$ 
for every $\theta\in\DD(\rho')$,
the complementary set of $D_n$ in $B_n$, contained in $C_n$, includes
all (but not only) vector fields with the same winding ratio as $\omega_n$.
\cqd

One can determine the non-constant part of each iterate $X_n=\omega_n+f_n$ by the
recurrence formula:
$$%\begin{equation}\label{recurrence formula f_n+1}
(\Ii-\Ee)\,f_{n+1}= \LL_n(\Ii-\Ee)\,f_{n} +
(\Ii-\Ee)\,\OO_n(\|f_n\|_{\rho'}^2),\quad n\geq0.
$$%\end{equation}
By induction, we obtain
\begin{equation}\label{formula f_n+1}
\begin{array}{rl}
(\Ii-\Ee)\,f_{n+1}= & \LL_n\circ\cdots\circ\LL_0\,(\Ii-\Ee)\, f_0 \\
& + \sum\limits_{j=1}^{n} \LL_{n}\circ\dots\circ\LL_j \,
(\Ii-\Ee)\,\OO_{j-1}\left(\|f_{j-1}\|_{\rho'}^2\right)  \\
& + (\Ii-\Ee)\,\OO_n\left(\|f_{n}\|_{\rho'}^2 \right).
\end{array}
\end{equation}

Note that the symbol $\OO_n(\|f_n\|^2_{\rho'})$ denotes the quadratic
remaining of $\RR_n$ as in Proposition \ref{proposition existence of
one-step renormalisation} which we write as $F_n^2$.

\begin{lemma}\label{recurrent formula f_n}
If $\omega_0\in DC(\beta)$, $\beta<1$,
there exists $\sigma,C,\tau>0$ such that, if
$\|f_0\|_{\rho'} \leq C^2<1$,
we have $\|f_{n}\|_{\rho'} \leq C \tilde A_{n}^{-\tau}<\zeta_n$, $n\geq1$.
\end{lemma}

\proof
We prove this lemma by induction.
From the formula $f_{n}=(\Ii-\Ee)f_{n} + \Pp_{n}\circ\Ee\,(f_{n}) +
(\Ii-\Pp_{n})\circ\Ee\,(f_{n})$, and referring to (\ref{P E f_n+1}) and Lemma \ref{lemma of the winding ratio cone},
$$
\|f_{n}\|_{\rho'} 
\leq  2 \|(\Ii-\Ee)f_{n}\|_{\rho'} + F_{n-1}^2,
$$
where
$F_{n-1}^2=\zeta_{n-1}^{-1}(\zeta_{n-1}-\|f_{n-1}\|_{\rho'})^{-1}
\|f_{n-1}\|_{\rho'}^2$ according to Proposition \ref{proposition
existence of one-step renormalisation}.
We write $\zeta_n=c'/(\alpha_n\alpha_{n+1})$ for some constant
$c'>0$.
Note that $C/\tilde A_n^\tau<\frac12 \zeta_n$ for a good choice of $C$
and $\tau$ independent of $n$, and by using the diophantine conditions.

Notice that $e^{-t} \leq (\nu/t)^{\nu}$ for any $t>0$ and
$\nu>0$.
Therefore, for $\tau>0$ and $0\leq j \leq n$, and the constants $c_1,c_2>0$
given in Proposition \ref{convergence on
stable mfld}:
\begin{equation}\label{exponential to quotient}
c_1 e^{-c_2 \Lambda_{j,n}} 
\leq 
c_1\sigma^\tau
\left[\frac{\tau(2+\beta)}{c_2}\right]^{\tau(2+\beta)}\frac{\tilde A_{j-1}^{\tau(1+\beta)}}{
\tilde A_{n+1}^\tau}  
\leq
\frac{\tilde A_{j-1}^{\tau(1+\beta)}}{7\tilde A_{n+1}^\tau},
\end{equation}
where we have assumed $\sigma^\tau\leq g(c_2)$, with a function
$g(y)= [y/(\tau(2+\beta))]^{\tau(2+\beta)}/(7c_1)$,
$y>0$, for the last inequality to hold.

If $n=1$, provided that $C<\frac12\zeta_0$ and
$2C^2(1+1/\zeta_0)<C/\tilde A_1^\tau<\frac12\zeta_1$, we have
$$
\|f_1\|_{\rho'} \leq
2C^2+\frac{C^2}{\zeta_0(\zeta_0-C)}
\leq
\frac {C}{\tilde A_1^\tau}<\frac12\zeta_1.
$$
Assuming the result is true for $n\geq1$, it remains to determine a bound of $\|f_{n+1}\|_{\rho'}$.
So, we look at each of the terms of the formula (\ref{formula f_n+1}):
\begin{equation}\label{induction for Ln..Lj a}
\begin{array}{rl}
\|(\Ii-\Ee)\,f_{n+1}\|_{\rho'}
\leq & 
\| \LL_n\circ\cdots\circ\LL_0\,(\Ii-\Ee)\, f_0\|_{\rho'}
\\
&
+ \sum\limits_{j=1}^{n}  \|\LL_{n}\circ\dots\circ\LL_j\,(\Ii-\Ee)\|.
F_{j-1}^2  
\\
& 
+ F_{n}^2.
\end{array}
\end{equation}
The first two are estimated by the use of Proposition \ref{convergence on
stable mfld} through (\ref{exponential to quotient}).
Therefore,
$$
\|\LL_n\circ\cdots\circ\LL_0\,(\Ii-\Ee)\,f_0\|_{\rho'}
\leq 
c_1 C e^{-c_2 \Lambda_{0,n}}
\leq 
\frac{C}{7\tilde A_{n+1}^\tau}.
$$
For the remaining terms we make use of the diophantine conditions in (\ref{diophantine conditions}) with
constants $K$, $\bar K$ and $\beta$: $\alpha_{n+1}\leq K\tilde
A_n^\beta$ and $\tilde A_{n+1}\leq \bar K\tilde A_n^{1+\beta}$.
So,
\begin{eqnarray*}
F_n^2
&\leq&
\frac{2C^2}{{c'}^2}\frac{\alpha_n^2\alpha_{n+1}^2}{\tilde A_n^{2\tau}}
\leq
\frac{2C^2K^2}{{c'}^2}\frac{\alpha_n^2\tilde A_{n}^{2\beta}}{\tilde A_n^{2\tau}}
\\
&\leq&
\frac{2C^2K^2\bar K^{(2\tau-2-2\beta)/(1+\beta)}}{{c'}^2}
\frac{1}{\tilde A_{n+1}^{(2\tau-2-2\beta)/(1+\beta)}}
\leq
\frac C{7\tilde A_{n+1}^\tau},
\end{eqnarray*}
for suitable $C>0$ and $\tau\geq(2\beta+2)/(1-\beta)$.
Now,
\begin{equation}\label{estimate Ln..Lj sum j}
\begin{array}{rcl}
\sum_{j=1}^n 
\|\LL_{n}\circ\dots\circ\LL_j \, (\Ii-\Ee)\|  F_{j-1}^2  
&\leq& 
\frac{2C^2K^2}{7{c'}^2}\sum_{j=1}^n
\frac{\tilde A_{j-1}^{\tau(1+\beta)+2\beta+2}}{\tilde
A_{n+1}^\tau\tilde A_{j-1}^{2\tau}}
\\
&\leq&
\frac{2C^2K^2}{7{c'}^2}\sum_{j=1}^n
\tilde A_{n+1}^{-\tau} \tilde A_{j-1}^{-t},
\end{array}
\end{equation}
if $\tau\geq(2\beta+2+t)/(1-\beta)$ for some $t>0$ and $\beta<1$.
Since $\tilde A_{j-1}^{-t}$ can be bounded by $\gamma^{-t(j-1)}$
times a constant, the sum in (\ref{estimate Ln..Lj sum j}) can be
bounded by some $M>0$.
Hence, (\ref{estimate Ln..Lj sum j}) is estimated by $C/(7\tilde
A_{n+1}^\tau)$ as long as $C$ is small.

Finally, if $C, \tau$ and $\sigma$ are chosen accordingly to the various conditions described
above, then
$$
\| f_{n+1}\|_{\rho'} \leq 
\frac{2(C/7 + C/7 + C/7) + C/7}{\tilde A_{n+1}^\tau}
=
\frac{C}{\tilde A_{n+1}^\tau}<\frac12\zeta_{n+1},
$$
as required.
\cqd

From the lower bound on $\tilde A_n$ in (\ref{bound on tilde A(n)}),
$\|f_n\|_{\rho'}$
decreases with $n$ at least geometrically like $\gamma^{-\tau n}$.
That completes
the proof of Theorem \ref{theorem generic frequency d=2}.

%%%%%%%%%%%%%%%%%%%%%%%%%%%%%%%%%%%%%%%%%%%%%%%%%%%%%%%%%%%%%%%%%%%%%%%%%%%%
\section*{Acknowledgements} 

I am indebted to Professor R. S. MacKay for 
the orientation and multiple suggestions given. 
I also wish to acknowledge Kostya Khanin and Nuno Rom\~ao for
useful discussions.
The author is supported by Funda\c c\~ao para a Ci\^encia e
a Tecnologia, under the research grant BD/11230/97.
%%%%%%%%%%%%%%%%%%%%%%%%%%%%%%%%%%%%%%%%%%%%%%%%%%%%%%%%%%%%%%%%%%%%%%%%%%%%

\bibliographystyle{plain} 
\addcontentsline{toc}{section}{References}
\bibliography{rfrncs}

\end{document}